\documentclass{amsart}
\usepackage{tikz-cd}
\usepackage{amsmath}
\usepackage{amsfonts}
\usepackage{mathtools}
\usepackage{amsmath,amssymb,amsthm,mathrsfs}
\usepackage[utf8]{inputenc}
\usepackage{enumitem}
\usepackage[pagebackref=true, colorlinks=true, linkcolor=blue]{hyperref}

\def\MR#1{}

\theoremstyle{plain} 
\newtheorem{theorem}{Theorem}
\numberwithin{theorem}{section}
\newtheorem*{theorem*}{Theorem}
\newtheorem{prop}[theorem]{Proposition}
\newtheorem{prop-def}[theorem]{Proposition-Definition}
\newtheorem{lemma}[theorem]{Lemma}
\newtheorem{coro}[theorem]{Corollary}
\newtheorem{conj}[theorem]{Conjecture}

\theoremstyle{definition}

\theoremstyle{remark} 
\newtheorem{remark}[theorem]{Remark}
\theoremstyle{definition} 
\newtheorem*{AIst}{AI statements}
\newtheorem*{ack}{Acknowledgment}

\newcommand{\wt}{\widetilde}

\newcommand{\be}{\begin{equation} }
\newcommand{\ee}{\end{equation} }

\usepackage{color}

\newcommand{\A}{\mathbb{A}}

\newcommand{\C}{\mathbb{C}} 
\newcommand{\D}{\mathbb{D}}

\renewcommand{\P}{\mathbb{P}}
\newcommand{\Q}{\mathbb{Q}}

\newcommand{\Z}{\mathbb{Z}}

\newcommand{\cE}{\mathcal{E}}

\newcommand{\cH}{\mathcal{H}}

\newcommand{\cM}{\mathcal{M}}
\newcommand{\cN}{\mathcal{N}}
\newcommand{\cO}{\mathscr{O}}

\newcommand{\cQ}{\mathcal{Q}}

\newcommand{\cT}{\mathcal{T}}

\newcommand{\sO}{\mathscr{O}}

\newcommand{\sA}{\mathscr{A}}

\newcommand{\shD}{\mathscr{D}}

\newcommand{\Spec}{\textup{Spec }}

\newcommand{\supp}{\textup{supp}}

\newcommand{\Ann}{{\textup{Ann}}}

\newcommand{\sRHom}{\mathcal R\mathcal Hom}

\title{The strong monodromy conjecture for hyperplane arrangements}
\author{Lei Wu}
\address{Lei Wu, School of Mathematical Sciences, Zhejiang
University, Hangzhou, China}
\email{leiwu23@zju.edu.cn}

\begin{document}

\maketitle

\begin{abstract}
    In this article, we prove the strong monodromy conjecture for complex hyperplane arrangements by proving a conjecture of Budur, Musta\c t\u a and Teitler that $-n/d$ is a root of the $b$-function of an irreducible essential and central hyperplane arrangement $f$ of degree $d$ on $\C^n$. 
\end{abstract}

\section{Introduction}\label{subsect:1.1}
We first recall the monodromy conjecture and the special case of hyperplane arrangements. Let $f\in \C[z_1,\dots,z_n]$ be a polynomial such that $f(0)=0$. Let $X = \A^n_\C$, and denote by $D=(f=0)$ the divisor of $f$, counting multiplicities. Take a log resolution $\mu\colon Y\to X$ of $D$. Define integers $a_i$ and $k_i$ by
\[ \mu^* D = \sum_{i\in S} a_i E_i, \quad K_{Y/X} = \sum_{i\in S} k_i E_i, \]
where the $E_i$ are irreducible divisors, and $K_{Y/X}$ is the relative canonical divisor of $\mu$.
The (local) \emph{topological zeta function} of $f$ is 
\[Z_f(s)\coloneqq \sum_{I\subseteq S} \chi(E^o_I\cap \mu^{-1}(0))\prod_{i\in I}\dfrac{1}{a_is+k_i+1},\]
where $E^o_I=\bigcap_{i\in I}E_i\setminus \bigcup_{j\in S \setminus I}E_j$ and $\chi$ is the topological Euler characteristic. Denef and Loeser \cite{DL} proved that $Z_f(s)$ is independent of the choice of log resolution.

The poles of $Z_f(s)$ are conjectured to have geometric and algebraic meaning. 
\begin{conj}[\cite{DL}, Topological Monodromy Conjecture]\label{conj:tmc}
If $c$ is a pole of the topological zeta function $Z_f(s)$, then:
\begin{itemize}
    \item $($Weak version$):$ $\exp ( 2\pi i c)$ is an eigenvalue of the monodromy action on the cohomology of the Milnor fiber.
    \item $($Strong version$):$ $c$ is a root of the Bernstein--Sato polynomial of $f$.
\end{itemize}
\end{conj}

It is well known that, by a classical result of Kashiwara and Malgrange, the strong version implies the weak version. The above conjecture is one of the most important conjectures in singularity theory. It remains wide open in general.

A testing ground for Conjecture \ref{conj:tmc} is the case of hyperplane arrangements, i.e., $f = \prod_{i} f_i$, where each $(f_i=0)$ defines a hyperplane in $\C^n$.
Budur, Musta\c{t}\u{a}, and Teitler \cite{BMT} studied Conjecture \ref{conj:tmc} in this case. They proved the weak version for hyperplane arrangements and also reduced the strong version to the following conjecture. We use $b_f(s)$ to denote the Bernstein--Sato polynomial (or $b$-function) of $f$. A hyperplane arrangement is \emph{central} if all hyperplanes pass through the origin $0$, and it is \emph{essential} if $\{0\}$ is a stratum of the intersection lattice. A central hyperplane arrangement is \emph{irreducible} if there is no linear change of coordinates on $\C^n$ such that $f$ can be written as a product of two non-constant polynomials in disjoint sets of variables. 

\begin{conj}[Budur, Musta\c{t}\u{a}, and Teitler]\label{conj:BMT}
If $f\in \C[z_1,\dots,z_n]$ is a polynomial of degree $d$ such that $(f=0)$ defines an irreducible, essential, central hyperplane arrangement, then $b_f(-\frac{n}{d})=0$.
\end{conj}
\begin{theorem}[Budur, Musta\c{t}\u{a}, and Teitler]\label{thm:BMT}
Conjecture \ref{conj:BMT} implies the strong version of Conjecture \ref{conj:tmc} for hyperplane arrangements.
\end{theorem}

Conjecture \ref{conj:BMT} is known for reduced $f$ when $n\le 3$ \cite{BSY}, for Weyl hyperplane arrangements \cite{BW}, for tame hyperplane arrangements \cite{Wal} (see also \cite[Thm.~1.3]{Bath2020} for the fact that the roots of $b_f$ in $[-1,0)$ are indeed combinatorial for tamed arrangements), for nonresonant hyperplane arrangements \cite{XYn/d}, and for $f$ with generic multiplicities \cite{SZ}. 
When the first version of this paper appeared, Davis and Yang \cite{YDn/d} claimed a proof of Conjecture \ref{conj:BMT} by developing Sabbah’s theory of multi-$V$-filtrations and using specialization along the diagonal.

In this paper, we prove Conjecture \ref{conj:BMT} in general, and hence Conjecture \ref{conj:tmc} for hyperplane arrangements, in a completely different (and much simpler) way. More precisely, we prove:
\begin{theorem}\label{thm:stmc}
Let $f$ be a polynomial of degree $d$ (not necessarily reduced) such that $(f=0)$ defines an irreducible, essential, central hyperplane arrangement. Then $-k/d$ is a root of $b_f(s)$ for $k=n,n+1,\dots,d$. In particular, Conjecture \ref{conj:BMT} holds.
\end{theorem}

The proof relies on the theory of algebraic relative (holonomic) $\shD$-modules. Specifically, we prove the commutativity of duality functors with proper direct images for algebraic relative $\shD$-modules (see Proposition \ref{prop:cmmdualpushfw}), which serves as an algebraic analogue of a fundamental result by Schapira and Schneiders \cite{SS94} for analytic relative $\shD$-modules. Another key ingredient is the Beilinson--Bernstein construction of $\shD$-module nearby cycles (see \S\ref{sect:3}). By combining this commutation property, the Beilinson--Bernstein nearby cycles, and resolution of singularities, we first provide an alternative proof of Kashiwara's rationality theorem for $b$-functions \cite{KasBf}, as well as Lichtin's refinement (Theorem \ref{thm:kashili}). The strategy underlying our proof of Theorem \ref{thm:kashili} establishes the framework for proving Theorem \ref{thm:stmc}. Furthermore, to better control the section ``$f^s$'' before and after passing to a resolution, we use the log canonical compactification of the hyperplane arrangement complement constructed in \cite{HKT2006}, together with the wonderful model, to construct an ample $\Q$-divisor whose round-up is the log canonical divisor $K_E+B$ and whose fractional part is supported on the boundary $B$. Here, $E$ denotes the wonderful model of the projectivized hyperplane arrangement in $\P^{n-1}$. We then apply the logarithmic Akizuki–Nakano vanishing theorem of \cite{AMPW2023} to show that the image of $f^s$ in the nearby cycles $\Psi_f^{-n/d}$ is nonzero. This approach also provides a general strategy for identifying roots of $b_f(s)$ in $(-1,0)$, whenever an explicit log resolution of $D$ is available.

\begin{AIst}
The AI-related Lemma 4.5 in the previous version of this paper is incorrect. GPT-6-Astra proposed the new Lemma \ref{lm:NTvanishing} as a replacement and suggested proving it using log canonical compactifications of hyperplane arrangement complements together with the logarithmic Akizuki–Nakano vanishing. Nevertheless, the remainder of the paper, including the overall strategy for proving Theorem \ref{thm:stmc}, remains valid and was developed independently of AI.

\end{AIst}

\begin{ack}
The author thanks Nero Budur for sharing his conjecture with him approximately five years ago. He is also grateful to Guillem Blanco Fernandez and Bradley Dirks for their valuable comments and suggestions on the early versions of this paper. 

\end{ack}

\section{Algebraic Relative $\shD$-modules}\label{sec:algrD}
We first recall the theory of algebraic relative $\shD$-modules. Let $X$ be a smooth algebraic variety over $\C$ and let $R$ be a Noetherian integral commutative $\C$-algebra such that $R_q$ is a regular local ring for every prime ideal $q\in \Spec R$. Set 
\[\sA^R=\sA^R_X\coloneqq \shD_X\otimes_\C R,\]
where $\shD_X$ is the sheaf of rings of algebraic differential operators. Similarly to $\shD_X$, $\sA^R$ is a coherent and Noetherian sheaf of rings.

If $\cN$ is a (left or right) $\sA^R$-module and $q\in \Spec R$, then since $R$ is central in $\sA^R$, 
\[\cN_q\coloneqq\cN\otimes_R R_q\]
becomes an $\sA^R_q \simeq \sA^{R_q}$-module, where $R_q$ is the localization of $R$ at $q$. We then define
\[\supp(\cN)=\supp^R(\cN)\coloneqq\{q\in \Spec R\mid \cN_q\neq 0\}\subseteq \Spec R\]
to be the support of $\cN$ as an $R$-module. 

Since $R$ is a $\C$-algebra, we have a natural inclusion $\shD_X\hookrightarrow \sA^R$, and thus $\sA^R$-modules are, in particular, $\shD_X$-modules. For a left $\sA^R$-module $\cM$, 
\[\cM^r\coloneqq \cM\otimes_{\cO_X}\omega_X\]
is the corresponding right $\shD_X$-module of $\cM$ under the side-change operation, where $\omega_X$ is the canonical line bundle. Since $R$ is central, $\cM^r$ is naturally a right $\sA^R$-module. Therefore, the side-change operation for $\shD$-modules extends naturally to $\sA^R$-modules.  
    
For a right $\sA^R$-module $\cN$ (or, more generally, a bounded complex of right $\sA^R$-modules), we define its dual by 
\[\D(\cN)\coloneqq \sRHom_{\sA^R}(\cN,\omega_X\otimes_{\cO_X}\sA^R)[\dim X];\]
similarly, for a left $\sA^R$-module $\cM$, we define
\[\D(\cM)\coloneqq \sRHom_{\sA^R}(\cM,\sA^R\otimes_{\cO_X}\omega_X^{-1})[\dim X].\]
It is clear that duality and the side-change operation commute, i.e.,
\[\D(\cM)\otimes_{\cO_X}\omega_X\simeq \D(\cM^r).\]

Let $g\colon Y\to X$ be a morphism between smooth algebraic varieties, and let $\cN^\bullet$ be a bounded complex of right $\sA^R_Y$-modules. Then the relative $\sA^R$-module direct image is defined as
\[g_+(\cN^\bullet)\coloneqq Rg_*(\cN^\bullet\otimes^L_{\sA_Y^R}g^*\sA^R_X)\simeq Rg_*(\cN^\bullet\otimes^L_{\shD_Y}g^*\shD_X),\]
where the isomorphism follows because $\sA_Y^R=\shD_Y\otimes_\C R$ and $g^*\sA^R_X\simeq g^*\shD_X\otimes_\C R$. This means that the relative direct image of $\cN^\bullet$ is precisely its direct image as a complex of $\shD$-modules.

\begin{prop}\label{prop:cmmdualpushfw}
    Let $g\colon Y\to X$ be a proper morphism between smooth algebraic varieties. Then we have a natural isomorphism of functors
    \[g_+\D\xrightarrow{\simeq} \D g_+\colon D^b_c(\sA^R_Y)^{\textup{op}} \longrightarrow D^b_c(\sA^R_X)^{\textup{op}},\]
    where $D^b_c(\sA^R_\bullet)^{\textup{op}}$ is the derived category of bounded complexes of right $\sA^R_\bullet$-modules with coherent cohomology sheaves. 
\end{prop}
\begin{proof}
    First, similar to the classical case (i.e., the case when $R=\C$), both $g_+$ and $\D$ preserve coherence (cf.\ \cite[Prop.~2.6.5 and Thm.~2.5.1]{HTT}).

    By definition, the trace map in \cite[\S 2.7.2]{HTT} provides the right-module version of the trace map:
    \[\textup{Tr}_g^{\textup{op}}\colon g_+(\omega_Y)[\dim Y]\longrightarrow \omega_X[\dim X].\]
    Similarly to the proof of \cite[Thm.~2.7.2]{HTT}, the trace map $\textup{Tr}_g^{\textup{op}}$ then induces a canonical isomorphism 
    \[g_+\D\xrightarrow{\simeq} \D g_+.\] 
\end{proof}
\begin{remark} 
    The above proposition can be seen as a simplified algebraic analogue of the fundamental result \cite[Thm.~5.15]{SS94} for analytic relative $\shD$-modules. 
\end{remark}

\section{Beilinson--Bernstein Nearby Cycles for $\shD$-modules}\label{sect:3}
In this section, we recall the Beilinson--Bernstein construction of $\shD$-module nearby cycles for regular functions on smooth algebraic varieties \cite[\S 4.2]{BBJats}. See also \cite[\S 3]{Gil} and \cite{BG, Wub}. We present the algebraic right-module version of the construction here for the application of this paper.

Let $X$ be a smooth algebraic variety over $\C$, and let $f$ be a non-invertible regular function on $X$. We consider the left $\sA^{\C[s]}$-module 
\[\cM\coloneqq j_*(\cO_U[s]f^s)=\cO_X[s,1/f]f^s\]
and $\cM^r=\cM\otimes_{\cO_X}\omega_X$, which is the corresponding right $\sA^{\C[s]}$-module under the side-change operation, where $j\colon U\coloneqq X\setminus D\hookrightarrow X$ is the open embedding and $D$ is the effective divisor defined by $(f=0)$. If $v$ is a global section of 
\[\omega_X(*D)\coloneqq\varinjlim_{k\to \infty} \omega_X(kD)\]
such that $v\cdot\sO_X|_U=\omega_U$, then we consider the coherent $\sA^{\C[s]}$-submodule 
\[\cN\coloneqq f^sv\cdot\shD_X[s]\subseteq \cM^r.\]
More generally, if $\cE\subseteq \omega_X(*D)$ is a coherent $\cO_X$-submodule such that $\cE|_U=\omega_U$, then we similarly obtain a coherent $\sA^{\C[s]}$-module
\[\cN(\cE)\coloneqq f^s\cE\cdot \shD_X[s]\subseteq\cM^r.\]
Then, for $l\in \Z$, we write
\[
\cN(lD)\coloneqq\cN(\cO_X(lD)\cdot v)=f^{s-l}v\cdot\shD_X[s]\subseteq \cM^r.
\]
It is clear that
\[\cN(lD)\hookrightarrow \cN((l+k)D)\]
for $k\in \Z_{\ge0}$.

The $b$-function of $f$ and $v$, denoted by $b_{f,v}$, is the monic generator of the annihilator ideal in $\C[s]$:
\[\Ann_{\C[s]}\left(\frac{\cN}{\cN(-D)}\right)\subseteq \C[s].\]
By substitution, we immediately see that for every $k\in \Z$,
\begin{equation}\label{eq:subk}
b_{f,v}(s+k)=b_{f,f^kv}(s).
\end{equation}
In the case when $X=\A^n_\C=\Spec \C[z_1,\dots,z_n]$ and $v=dz_1\wedge dz_2\wedge \dots \wedge dz_n$, since the action of $\C[s]$ and the side-change operation commute, $b_{f,v}$ is exactly the Bernstein--Sato polynomial $b_f$ of $f\in \C[z_1,\dots,z_n]$.

\begin{lemma}\label{lm:z=bf}
    With the notation as above, we have 
    \[Z(b_{f,v})=\supp\left(\frac{\cN}{\cN(-D)}\right),\]
    where $Z(b_{f,v})\subseteq \Spec\C[s]$ is the algebraic set defined by $b_{f,v}$, i.e., the set of all roots of $b_{f,v}$.
\end{lemma}
\begin{proof}
    Take a finite affine open covering $X=\bigcup_i V_i$ such that $\omega_{V_i}\simeq \cO_{V_i}\cdot e_i$ and $v|_{V_i}=h_i\cdot e_i$ for some rational function $h_i$ regular on $V_i\setminus D$. By \cite[R\'esultat 2]{Mai}, 
    \[\frac{\shD_{V_i}[s]h_i\cdot f^s}{\shD_{V_i}[s]h_i\cdot f^{s+1}}\] 
    is algebraic relative holonomic (cf.\ \cite[Def.~3.2.3]{BVWZ}). 
    By \cite[Lem.~3.4.1]{BVWZ}, we know 
    \[Z(b_{f,h_i})=\supp\left(\frac{\shD_{V_i}[s]h_i\cdot f^s}{\shD_{V_i}[s]h_i\cdot f^{s+1}}\right),\]
    where $b_{f,h_i}$ is the monic generator of $\Ann_{\C[s]}\left(\frac{\shD_{V_i}[s]h_i\cdot f^s}{\shD_{V_i}[s]h_i\cdot f^{s+1}}\right)$. Applying the side-change operation, we then obtain 
    \[Z(b_{f|_{V_i},v|_{V_i}})=Z\left(\Ann_{\C[s]}\left(\left.\frac{\cN}{\cN(-D)}\right|_{V_i}\right)\right)=\supp\left(\left.\frac{\cN}{\cN(-D)}\right|_{V_i}\right).\]
    By definition, we know
    \[\Ann_{\C[s]}\left(\frac{\cN}{\cN(-D)}\right)=\bigcap_i \Ann_{\C[s]}\left(\left.\frac{\cN}{\cN(-D)}\right|_{V_i}\right)\]
    and 
    \[\supp\left(\frac{\cN}{\cN(-D)}\right)=\bigcup_i \supp\left(\left.\frac{\cN}{\cN(-D)}\right|_{V_i}\right).\]
    Since $\C[s]$ is a PID, we thus obtain
    \[Z(b_{f,v})=\supp\left(\frac{\cN}{\cN(-D)}\right).\]
\end{proof}

It is well known that the local $b$-function $b_{f|_{V_i},v|_{V_i}}$ is always nonzero and non-constant if $f|_{V_i}$ is not invertible. The above proof also shows that $b_{f,v}$ is the least common multiple of all the $b_{f|_{V_i},v|_{V_i}}$; in particular, the global $b$-function $b_{f,v}$ is non-constant.  

Fixing an arbitrary $\alpha\in \C\simeq \Spec \C[s]$, we use $\mathfrak{m}=\mathfrak{m}_\alpha\subseteq \C[s]$ to denote the maximal ideal corresponding to $\alpha$. 
We define
\[j_!(\cM^r_{\mathfrak{m}}|_U)\coloneqq \D j_*\big(\D(\cM^r_{\mathfrak{m}}|_U)\big)\]
and 
\[j_!(\cM^r_\eta |_U)\coloneqq \D j_*\big(\D(\cM^r_\eta|_U)\big),\]
where $\eta\in \Spec \C[s]$ is the generic point. 

Since taking duals and the side-change operation commute, by \cite[Thm.~2.10 and Prop.~2.12(2)]{Wub}, we have the following right-module version of a fundamental result of Beilinson--Bernstein; see also \cite[\S 4.2]{BBJats}, \cite[\S 3.6]{Gil}, and \cite{BG}. 

\begin{theorem}[Beilinson--Bernstein]\label{thm:BBright}
    With the notation and assumptions as above, we have:
    \begin{enumerate}
        \item $j_!(\cM^r_{\mathfrak{m}}|_U)\stackrel{\textup{q.i.}}{\simeq}\cH^0\big(j_!(\cM^r_{\mathfrak{m}}|_U)\big)$, i.e., $j_!(\cM^r_{\mathfrak{m}}|_U)$ is a sheaf;
        \item the natural morphism $j_!(\cM^r_{\mathfrak{m}}|_U)\hookrightarrow j_*(\cM^r_{\mathfrak{m}}|_U)$ is injective;
        \item $j_*(\cM^r_{\mathfrak{m}}|_U)=\cM^r_{\mathfrak{m}}=\cN(lD)_{\mathfrak{m}} $ for $l\gg 0$;
        \item $j_!(\cM^r_{\mathfrak{m}}|_U)=\cN(-lD)_{\mathfrak{m}}$ for $l\gg 0$;
        \item forgetting the $\C[s]$-module structure, $\Psi_{f}^\alpha\coloneqq\frac{j_*(\cM^r_{\mathfrak{m}}|_U)}{j_!(\cM^r_{\mathfrak{m}}|_U)}$ is a regular holonomic $\shD_X$-module;
        \item $j_!(\cM^r_\eta|_U)= j_*(\cM^r_\eta|_U)$. 
    \end{enumerate}
\end{theorem}

The sheaf $\Psi_{f}^\alpha$ is the right-module version of the $\alpha$-nearby cycles of $f$, which corresponds to the topological nearby cycles of the constant sheaf $\C_X$ along $f$ under the Riemann--Hilbert correspondence (see \cite[\S 4.2]{BBJats} and \cite[Thm.~1.1]{Wub}).

By Theorem \ref{thm:BBright}, $\frac{\cN(kD)_{\mathfrak{m}}}{\cN((k-1)D)_{\mathfrak{m}}}$ is always a subquotient of $\Psi^\alpha_f$ for every $k\in \Z$. Applying Eq.~\eqref{eq:subk} and Lemma \ref{lm:z=bf} together gives:

\begin{coro}\label{cor:bfest}
    For $\alpha\in\C$, the condition that $b_{f,v}(\alpha+k)=0$ for some $k\in \Z$ is equivalent to $\Psi^\alpha_f\neq 0$.
\end{coro}

\section{Proof of the main result}
We fix a non-invertible regular function $f$ on $X=\A^n_\C$ and a log resolution $\mu\colon Y\to X$, along with the notation from \S\ref{subsect:1.1}. Kashiwara \cite{KasBf} proved that the roots of $b_f$ are negative rational numbers. Lichtin refined Kashiwara's result into the following form:
\begin{theorem}\cite[Thm.~5]{Lich}\label{thm:kashili}
    All the roots of $b_f$ are of the form 
    \[-\frac{k_i+l+1}{a_i}\]
    for $i\in S$ and $l\in \Z_{\ge0}$.
\end{theorem}
As a warm-up, we first give an alternative proof of the above theorem using the nearby cycles constructed in \S\ref{sect:3}.
\begin{proof}
    We set $\widetilde{f}=\mu^*f$, $v=dz_1\wedge dz_2\wedge \dots \wedge dz_n$, and $\widetilde{v}=\mu^*v\in \mu^*\omega_X\hookrightarrow\omega_Y$. We consider 
    \[\widetilde{\cN}\coloneqq \widetilde{f}^s\widetilde{v}\cdot \shD_Y[s]\subseteq \widetilde{j}_*(\cM^r|_U),\]
    where $\widetilde{j}\colon U\hookrightarrow Y$ is the open embedding. Then we have a commutative diagram 
    \[
    \begin{tikzcd}
        U\arrow[r,"\widetilde{j}"]\arrow[d,"\textup{id}"]& Y \arrow[d,"\mu"]\\
        U\arrow[r,"j"]& X
    \end{tikzcd}
    \]
    Since both $j$ and $\widetilde{j}$ are affine open embeddings, we have $j_+=j_*$ and $\widetilde{j}_+=\widetilde{j}_*$. By functoriality, the above diagram yields that for every $\alpha\in \C$,
    \[\mu_+(\widetilde{j}_*(\cM^r_{\mathfrak{m}}|_U))\simeq j_*(\cM^r_{\mathfrak{m}}|_U),\]
    where $\mathfrak{m}=\mathfrak{m}_\alpha$ is the maximal ideal corresponding to $\alpha$. Similarly, by Proposition \ref{prop:cmmdualpushfw}, we also know
    \begin{equation}\label{eq:birapushforminex}
    \mu_+(\widetilde{j}_!(\cM^r_{\mathfrak{m}}|_U))\simeq j_!(\cM^r_{\mathfrak{m}}|_U).
    \end{equation}
    By Theorem \ref{thm:BBright}(2), we have a short exact sequence
    \[0\to \widetilde{j}_!(\cM^r_{\mathfrak{m}}|_U)\longrightarrow\widetilde{j}_*(\cM^r_{\mathfrak{m}}|_U)\longrightarrow \Psi_{\widetilde{f}}^\alpha\to 0.\]
    Applying $\mu_+$ yields 
    \begin{equation}\label{eq:pushfnearby}
    \mu_+(\Psi_{\widetilde{f}}^\alpha)\simeq \Psi_{f}^\alpha. 
    \end{equation}

    We take a coordinate system $(y_1,\dots,y_n)$ on an open neighborhood $V\subseteq Y$ locally around some $p\in E$ such that 
    \[\widetilde{f}|_V=u\cdot\prod_{i\in S_p\subseteq S}y_i^{a_i}\quad \textup{and} \quad \widetilde{v}|_V=u'\prod_{i\in S_p}y_i^{k_i}dy_1\wedge dy_2\wedge \dots \wedge dy_n\]
    with $u,u'$ being invertible regular functions on $V$. Then a local calculation shows that the roots of $b_{\widetilde{f},\widetilde{v}}$ are of the form 
    \[-\frac{k_i+l+1}{a_i}\]
    for $i\in S$ and $l\in \Z_{\ge0}$. Corollary \ref{cor:bfest} and Equation \eqref{eq:pushfnearby} together imply that the roots of $b_f=b_{f,v}$  
    are of the form 
    \[-\frac{k_i+l+1}{a_i}\]
    for $i\in S$ and $l\in \Z$. 
    
    Now, we need to exclude all $l < 0$. By Theorem \ref{thm:BBright}(4), we have 
    \[\widetilde{j}_!(\cM^r_{\mathfrak{m}}|_U)=\widetilde{\cN}(-k\mu^*D)_{\mathfrak{m}}\]
    for $k\gg 0$. Suppose 
    \[\alpha\notin
\bigcup_{i\in S}
\left\{-\frac{k_i+l+1}{a_i}:l\ge0\right\}.\]
    Then $\alpha$ is not a root of $b_{\widetilde{f},\widetilde{v}}(s+k)$ for any $k\ge 0$. By Equation \eqref{eq:subk} and Lemma \ref{lm:z=bf}, we know 
    \begin{equation}\label{eq:spj!}
    \widetilde{j}_!(\cM^r_{\mathfrak{m}}|_U)=\widetilde{\cN}_{\mathfrak{m}}.
    \end{equation}
    Now, we consider the global section $\widetilde{f}^s\widetilde{v}\otimes 1$ of the sheaf $\widetilde{\cN}_{\mathfrak{m}}\otimes_{\shD_Y} \mu^*\shD_X$ (and the sheaf $\widetilde{\cN}_{\mathfrak{m}}\otimes_{\cO_Y} \mu^*\shD_X$),
    which induces morphisms in the derived category of $\mu^{-1}\cO_X$-modules
    \[\cO_Y\to\widetilde{\cN}_{\mathfrak{m}}\otimes_{\cO_Y}^L\mu^*\shD_X\to\widetilde{\cN}_{\mathfrak{m}}\otimes_{\shD_Y}^L\mu^*\shD_X\to \widetilde{j}_*(\cM^r_{\mathfrak{m}}|_U)\otimes^L_{\shD_Y} \mu^*\shD_X,\]
    since $\shD_Y$ is flat over $\cO_Y$. Applying $R\mu_*$, we further have the induced morphisms in the derived category of $\cO_X$-modules
    \[\cO_X\to R\mu_*\cO_Y\to R\mu_*\big(\widetilde{j}_*(\cM^r_{\mathfrak{m}}|_U)\otimes^L_{\shD_Y} \mu^*\shD_X\big)\simeq j_*(\cM^r_{\mathfrak{m}}|_U).\]
    Then $1\in \cO_X(X)$ yields a global section $w$ of $j_*(\cM^r_{\mathfrak{m}}|_U)$.
    Clearly, $w|_U=f^sv|_U$, and thus $w-f^sv$ is a section of $j_*(\cM^r_{\mathfrak{m}}|_U)$ supported on $D$. Since $j_*(\cM^r_{\mathfrak{m}}|_U)$ has no nonzero sections supported on $D$ by construction, we know $w=f^sv$. This means $w=f^sv$ is the image of $\tilde f^s\tilde v$ under $\mu_+$. 
    
    If such an $\alpha$ is a root of $b_{f,v}$, then we have 
    \[\frac{\cN_{\mathfrak{m}}}{\cN(-D)_{\mathfrak{m}}}\neq 0\]
    by Lemma \ref{lm:z=bf}, and thus by Theorem \ref{thm:BBright}(4),
    \[\frac{\cN_{\mathfrak{m}}}{j_!(\cM^r_{\mathfrak{m}}|_U)}\neq 0.\]
    Since $\cN_{\mathfrak{m}}$ is generated by $f^sv$, we know in particular that $f^sv\notin j_!(\cM^r_{\mathfrak{m}}|_U)$. However, by Equation \eqref{eq:spj!}, we have a factorization 
    \[\begin{tikzcd}
        \cO_X\arrow[rr]\arrow[rd]&& j_*(\cM^r_{\mathfrak{m}}|_U)\\
        & j_!(\cM^r_{\mathfrak{m}}|_U)\arrow[ru, hook].&
    \end{tikzcd}\]
    Then $w=f^sv$ is also a section of $j_!(\cM^r_{\mathfrak{m}}|_U)$, which is a contradiction. Consequently, $\alpha= -\frac{k_i+l+1}{a_i}$ is not a root of $b_f=b_{f,v}$ for all $l<0$. 
\end{proof}

\begin{coro}\label{cor:j!-10}
    For every $\alpha\in \Q\cap[-1,0)$, we have 
    \[j_!(\cM^r_{\mathfrak{m}}|_U)=\cN(-D)_{\mathfrak{m}}=f^{s+1}v\cdot\shD_X[s]_{\mathfrak{m}},\]
    where $\mathfrak{m}\subseteq \C[s]$ is the maximal ideal corresponding to $\alpha$.
\end{coro}
\begin{proof}
    By Theorem \ref{thm:BBright}(4), it suffices to prove 
    \[\frac{\cN(-kD)_{\mathfrak{m}}}{\cN(-(k+1)D)_{\mathfrak{m}}}=0\]
    for all $k\in \Z_{>0}$. If $\frac{\cN(-kD)_{\mathfrak{m}}}{\cN(-(k+1)D)_{\mathfrak{m}}}\neq 0$ for some $k>0$, then by Lemma \ref{lm:z=bf}, Equation \eqref{eq:subk}, and substitution, we conclude 
    \[b_f(\alpha+k)=0.\]
    Since $\alpha\in \Q\cap[-1,0)$ and $k\ge 1$, $\alpha+k\ge0$ is a root of $b_f(s)$, which contradicts Theorem \ref{thm:kashili}.
\end{proof}

Now, we prove Theorem \ref{thm:stmc}. Assume $f\in \C[z_1,\dots,z_n]$ is a polynomial of degree $d$ such that $(f=0)$ defines an irreducible, essential, central hyperplane arrangement $A$. Let $L(A)$ be the intersection lattice of $A$, and let $\widetilde{L}(A)$ denote the set of all dense edges in $L(A)$ (see, for instance, \cite{BMT} for definitions).
It is clear that $d>n$. Take $\alpha=-k/d$ for some $k\in \{n,n+1,\dots,d-1\}$. We may further assume 
\[f=\prod_{i}f_i^{d_i}\]
such that each $f_i$ is a linear polynomial. For each $W\in \widetilde{L}(A)$, we define 
\[d_W=\sum_{f_i|_W\equiv 0}d_i \quad \textup{and} \quad f_W=\prod_{f_i|_W\equiv 0}f_i^{d_i}.\]
We also set $r_W$ to be the codimension of $W\in L(A)$. For example, $d_{\{0\}}=d$ and $r_{\{0\}}=n$. We define 
\[\Lambda_\alpha=\left\{W\in \widetilde{L} (A)\mathrel{\Big|} \alpha=-k/d_W \textup{ for some } k\in \{r_W,r_W+1,\dots,d_W-1\}\right\}.\]
Since $\{0\}\in \Lambda_\alpha$, $\Lambda_\alpha$ is nonempty. We then take a $W_0\in \Lambda_{\alpha}$ such that 
\[r_{W_0}=\min\{r_W\mid W\in \Lambda_\alpha\}.\]
Since the global $b$-function is always a multiple of the local one, by discarding all hyperplanes $(f_i=0)$ that do not contain $W_0$, we can replace $f$ by $f_{W_0}$ and $X$ by $X/W_0$, keeping the same notation. Therefore, we know that if $W$ is a dense edge and $W\neq \{0\}$, then $\alpha\neq -k/d_W$ for any $k\in\{r_W,r_W+1,\dots,d_W-1\}$; i.e., $\Lambda_\alpha=\{\{0\}\}$. 

Let $\mu\colon Y\to X$ be the log resolution obtained by successively blowing up the (strict transforms of the) unions of the dense edges of dimension $m$ for $m=0,1,\dots,n-2$. By \cite[Thm.~3.1]{STV}, $\mu$ is indeed a log resolution of $A$. Then 
\[K_{Y/X}=\sum_{W\in \widetilde{L}(A)}(r_W-1)E_W \quad \textup{and} \quad \mu^*D=\sum_{W\in \widetilde{L}(A)}d_WE_W.\]
Set $\widetilde{f}=\mu^*(f)$, $v=dz_1\wedge\dots\wedge dz_n$, and $\widetilde{v}=\mu^*v$. As in the proof of Theorem \ref{thm:kashili}, we have 
\[\widetilde{\cN}=\widetilde{f}^s\widetilde{v}\cdot \shD_Y[s]\subseteq \widetilde{j}_*(\cM^r|_U).\]
Then, in a small open neighborhood $V=V_p\subseteq Y$ of a point $p\in \bigcup_{W\in \widetilde{L}(A)}E_W$, we have 
\begin{equation}\label{eq:localwtfhyp}
\widetilde{\cN}|_V=u^su'\prod_{W\in \widetilde{L}_p(A)} y_{W}^{d_Ws+r_W-1}dy\cdot \shD_V[s],
\end{equation}
where $u,u'$ are invertible regular functions on $V$, $y_W$ is the local coordinate defining $E_W$, $dy$ is the local volume form, and $\widetilde{L}_p(A)\coloneqq\{W\in \widetilde{L}(A)\mid E_W\cap V\neq \emptyset\}$. Since $u$ and $u'$ are invertible, the gauge-transformation given by the conjugation-action by $u'u^{s}$
induces an isomorphism of sheaves of rings 
\[\shD_V[s]\simeq u^{-s}(u')^{-1}\cdot \shD_V[s]\cdot u'u^{s},\]
which further induces an isomorphism
\[\widetilde{\cN}|_V\simeq\prod_{W\in \widetilde{L}_p(A)} y_{W}^{d_Ws+r_W-1}dy\cdot \shD_V[s]\] 
as $\shD_V[s]$-modules.
Therefore, we may safely assume $u=u'=1$ for simplicity.

For an arbitrary rational number $\beta\in \Q$, we use $\mathfrak{n}\subseteq \C[s]$ to denote the maximal ideal corresponding to $\beta$. Then, by Theorem \ref{thm:BBright}, we have  
\[\widetilde{j}_!(\cM^r_{\mathfrak{n}}|_U)=\widetilde{f}^{s+k}\widetilde{v}\cdot \shD_Y[s]_{\mathfrak{n}}\subseteq \widetilde{j}_*(\cM^r_{\mathfrak{n}}|_U)\]
for $k\gg 0$. Since $\alpha\neq -k/d_W$ for any $k\in\{r_W,r_W+1,\dots,d_W-1\}$ and for $W\not=\{0\}$, the local expression \eqref{eq:localwtfhyp} and Theorem \ref{thm:BBright}(4) yield 
\[\widetilde{j}_!(\cM^r_{\mathfrak{m}}|_U)=\widetilde{f}^s\cO_Y(-(k-n+1)E_{\{0\}})\widetilde{v}\cdot\shD_Y[s]_{\mathfrak{m}},\]
where $\mathfrak{m}$ is the maximal ideal corresponding to $\alpha=-k/d$. Therefore, we have a short exact sequence
\begin{equation}\label{eq:sesupstairnearby}
0\to \cQ_0\coloneqq\frac{\widetilde{\cN}_{\mathfrak{m}}}{\widetilde{j}_!(\cM^r_{\mathfrak{m}}|_U)}\longrightarrow \Psi_{\widetilde{f}}^\alpha\longrightarrow \frac{\widetilde{j}_*(\cM^r_{\mathfrak{m}}|_U)}{\widetilde{\cN}_{\mathfrak{m}}}\to 0.
\end{equation} 

By Theorem \ref{thm:BBright}(3) and \eqref{eq:localwtfhyp}, locally around $p\in V$, the quotient $\frac{\widetilde{j}_*(\cM^r_{\mathfrak{m}}|_U)}{\widetilde{\cN}_{\mathfrak{m}}}$ has a finite filtration such that each successive quotient is of the form 
\[\frac{y_{W'}^{-1}\prod_{W\in \widetilde{L}_p(A)} y_{W}^{d_Ws+r_W-1-c_W}dy\cdot\shD_V[s]_{\mathfrak{m}}}{\prod_{W\in \widetilde{L}_p(A)} y_{W}^{d_Ws+r_W-1-c_W}dy\cdot\shD_V[s]_{\mathfrak{m}}}\]
with $ c_W< r_{W}-1$ for every $W\in \widetilde{L}_p(A)$, where $W'$ runs over the set $\widetilde{L}(A)$. We only require $c_W< r_{W}-1$ because $\alpha<0$. By construction, we know that a successive quotient with respect to $W'$ is nonzero if and only if $-\frac{r_{W'}-1-c_{W'}}{d_{W'}}=\alpha$. In particular, all the successive quotients with respect to the edge $\{0\}$ are zero. Discarding all zero quotients, we isolate a nonzero successive quotient $\cQ_{W'}$ (depending only on $W'$) such that 
\begin{equation}\label{eq:localquotient}
\cQ_{W'}\simeq\frac{y_{W'}^{-1}\prod_{W\in \widetilde{L}_p(A)} y_{W}^{d_Ws+r_W-1-c_W}dy\cdot\shD_V[s]}{\prod_{W\in \widetilde{L}_p(A)} y_{W}^{d_Ws+r_W-1-c_W}dy\cdot\shD_V[s]}
\end{equation}
for $W'\in \Gamma_\alpha\coloneqq\left\{W\in \widetilde{L}(A)\mathrel{\big|} -\frac{r_{W}-1-c_{W}}{d_{W}}=\alpha \textup{ for some integer } 0 \le c_W<r_W-1\right\}$; necessarily we have $c_{W'} = d_{W'}\alpha +r_{W'} - 1$, and for the other $c_W$ one can make suitable choices. 
Similarly, for $p\in E_{\{0\}}$ we have 
\begin{equation}\label{eq:localisoQ}
\cQ_0|_V\simeq \frac{y^{k-n}_{\{0\}}\prod_{\substack{W'\in \widetilde{L}(A) \\ W'\neq \{0\}}}y_{W'}^{a_{W'}}\prod_{W\in \widetilde{L}_p(A)} y_{W}^{d_Ws+r_W-1}dy\cdot\shD_V[s]}{y^{k-n+1}_{\{0\}}\prod_{\substack{W'\in \widetilde{L}(A) \\ W'\neq \{0\}}}y_{W'}^{a_{W'}}\prod_{W\in \widetilde{L}_p(A)} y_{W}^{d_Ws+r_W-1}dy\cdot\shD_V[s]},
\end{equation}
where we can put $a_{W'}$ to be any nonnegative ingeter for every $W'\in \widetilde{L}(A)\setminus\{\{0\}\}$ since $\Lambda_\alpha=\{\{0\}\}$ by assumption. Different choice of $a_{W'}$ means we pick different local generators for $\cQ_0$. Since $\cQ_0|_V$ is already supported at ${\mathfrak{m}}$, it is not necessary to take localization.  

\begin{lemma}\label{lm:pushlocalization}
    With the notation and assumptions as above, we have 
    \[\cH^i\mu_+\left(\frac{\widetilde{j}_*(\cM^r_{\mathfrak{m}}|_U)}{\widetilde{\cN}_{\mathfrak{m}}}\right)=0\]
    for all $i<0$.
\end{lemma}
\begin{proof}
    Locally around $p$, we consider $\cQ_{W'}$ for every $W'\in \Gamma_\alpha$. By Kashiwara's equivalence and Equation \eqref{eq:localquotient}, $\cQ_{W'}$ corresponds to a (left) regular holonomic $\shD_{E_{W'}}$-module generated by $\bar e$ subject to the relations
    \[(y_W\partial_{y_W})\cdot \bar e=(d_W\cdot \alpha+r_W-1-c_W)\bar e\]
    for $W\neq W'$. In particular, if $\{0\}\in \widetilde{L}_p(A)$, then 
    \[(y_{\{0\}}\partial_{y_{\{0\}}})\cdot \bar e=n-1-k-c_{\{0\}}<0.\]
    Therefore, $\cQ_{W'}$ is an $\cO_Y(*E_{\{0\}})=\varinjlim \cO_Y(kE_{\{0\}})$-module. By construction, we then know $\frac{\widetilde{j}_*(\cM^r_{\mathfrak{m}}|_U)}{\widetilde{\cN}_{\mathfrak{m}}}$ is also an $\cO_Y(*E_{\{0\}})$-module. Since $\cQ_0$ is supported on $E_{\{0\}}$, using \eqref{eq:sesupstairnearby} we conclude 
    \[\Psi_{\widetilde{f}}^\alpha\otimes_{\cO_Y} \cO_Y(*E_{\{0\}})\simeq \frac{\widetilde{j}_*(\cM^r_{\mathfrak{m}}|_U)}{\widetilde{\cN}_{\mathfrak{m}}}\]
    and 
    \[\cQ_0\simeq R\Gamma_{[E_{\{0\}}]}(\Psi_{\widetilde{f}}^\alpha).\]
    By \cite[Ch.~II.2.5.28]{Bj} and \eqref{eq:pushfnearby}, we have a distinguished triangle
    \[R\Gamma_{[\{0\}]}(\Psi_f^\alpha)\longrightarrow \Psi_f^\alpha\longrightarrow \mu_+\left(\frac{\widetilde{j}_*(\cM^r_{\mathfrak{m}}|_U)}{\widetilde{\cN}_{\mathfrak{m}}}\right)\xrightarrow{+1}.\]
    Taking the associated long exact sequence, the required statement thus follows. 
\end{proof}

For simplicity, we write $E=E_{\{0\}}$, $E^o=E_{\{0\}}\setminus \bigcup_{W\neq \{0\}} E_W$ and $B=E\setminus E^o$.
\begin{lemma}\label{lm:NTvanishing}
    With notation as above, we have 
    \begin{equation}\label{eq:allvanish}
    H^j\big(E, \Omega^i_E(\log B)\otimes \omega_E(B)^{-1}\big)=0,\qquad i+j<n-1=\dim E.
    \end{equation}
    In particular, 
    \[H^{r-1}\big(E, \bigwedge^{r} \cT_{E}(-\log B)\big)
= 0, \qquad 1 \le r \le n-1,\]
where $\cT_{E}(-\log B)$  is the $\cO_E$-dual of $\Omega^1_E(\log B)$, the sheaf of log one-forms. 
\end{lemma}
\begin{proof}
Let $(S,B_S)$ be the log canonical model of
$U=E^o$, the complement of the projectivized hyperplane arrangement,
constructed in \cite[\S 2]{HKT2006}. The hypotheses that $A$ is central,
essential, and irreducible ensure that this construction applies.
By \cite[Thm.~2.2]{HKT2006}, $S$ is a normal projective variety,
the pair $(S,B_S)$ has toric singularities, and
\[
 A_S\coloneqq K_S+B_S
\]
is a very ample Cartier divisor. As explained in \cite[\S 4]{Tev07},
there is a log crepant birational morphism $\gamma\colon E\to S$
restricting to the identity on $U$, with
\begin{equation}\label{eq:crepant}
 L\coloneqq K_E+B\sim\gamma^*A_S,
 \qquad
 \bigl(\gamma^{-1}(B_S)\bigr)_{\mathrm{red}}=B.
\end{equation}

We now construct an ample $\Q$-divisor $A'$ such that
$\lceil A'\rceil=L$ and the fractional part of $A'$ is supported
on $B$. In the construction of \cite[\S 2]{HKT2006}, $S$ is embedded
in a normal projective toric variety $T$, and $B_S$ is the reduced
boundary induced by the toric boundary $B_T$ of $T$. Choose an ample
effective Cartier divisor $D_T$ with support exactly $B_T$.
Its restriction
\[
 D_S\coloneqq D_T|_S
\]
is an ample effective Cartier divisor with support exactly $B_S$.
Consequently, $\gamma^*D_S$ is an effective, nef, and big Cartier
divisor with support exactly $B$.

Choose an ample integral divisor $H_E$ supported on $B$.
Such a representative exists because $\operatorname{Pic}(E)$ is
generated by the pullback of the hyperplane class and the exceptional
divisor classes: the hyperplane class can be represented by an
arrangement hyperplane, whose total transform is supported on $B$.
Since $\gamma^*D_S$ has positive coefficients along every irreducible
component of $B$, we may choose an integer $M\gg 0$ such that
\[
 F_0\coloneqq M\gamma^*D_S-H_E
\]
is effective with support exactly $B$.

By openness of the ample cone, we may choose a sufficiently small
positive rational number $\varepsilon$ such that
$A_S-\varepsilon M D_S$ remains ample and every nonzero coefficient
of
\[
 \Delta\coloneqq\varepsilon F_0
\]
lies in $(0,1)$. Set $A'\coloneqq L-\Delta$. Then
\[
 A'\sim_{\Q}
 \gamma^*(A_S-\varepsilon M D_S)+\varepsilon H_E.
\]
The right-hand side is the sum of a nef $\Q$-divisor and an ample
$\Q$-divisor, so $A'$ is ample. Moreover, since $L$ is integral,
we have the following identities and inequality of divisors:
\[
 \lceil A'\rceil=L,
 \qquad
 \lceil A'\rceil-A'=\Delta\leq B.
\]
In particular, the fractional part of $A'$ is supported on $B$.

Set $m\coloneqq n-1=\dim E$. Applying the logarithmic
Akizuki--Nakano vanishing theorem
\cite[Thm.~1.3]{AMPW2023} to the ample $\Q$-divisor $A'$,
with both the logarithmic boundary and the integral divisor
subtracted in that theorem equal to $B$, gives
\[
 H^{m-j}\bigl(
 E,\Omega_E^{m-a}(\log B)(L-B)
 \bigr)=0
 \qquad\text{whenever}\qquad
 (m-j)+(m-a)>m.
\]
By Serre duality and the logarithmic wedge pairing, the dual of
this cohomology group is naturally identified with
\[
 H^j\bigl(E,\Omega_E^a(\log B)(-L)\bigr).
\]
Since $(m-j)+(m-a)>m$ is equivalent to $a+j<m$, we conclude that
\[
 H^j\bigl(E,\Omega_E^a(\log B)(-L)\bigr)=0
 \qquad\text{for }a+j<m,
\]
which proves \eqref{eq:allvanish}.
\end{proof}

By Kashiwara's equivalence and \eqref{eq:localisoQ}, $\cQ_0$
corresponds to a regular holonomic right $\shD_E$-module $\cE$
generated by the global section
\[
 e\coloneqq
 [\mu^*(f_1^{k-n})\cdot\tilde f^s\tilde v].
\]
Locally, this section satisfies the relations
\begin{equation}\label{eq:raction}
 e\cdot(-\partial_{y_W}y_W)
 =\bigl(d_W\cdot\alpha+r_W-1+a_W\bigr)e,
 \qquad W\neq\{0\},
\end{equation}
where each $a_W$ is a nonnegative integer. 
More precisely, $a_W$ is the multiplicity of $\mu^*(f_1^{k-n})$ along $y_W$.
The assumption on $\Lambda_\alpha$ implies that
$d_W\cdot\alpha+r_W-1$ is either nonintegral or a nonnegative
integer. Since $a_W$ is a nonnegative integer, the same holds
for the residue $d_W\cdot\alpha+r_W-1+a_W$.

The free $\cO_E$-module
\[
 \bar\cE\coloneqq\cO_E\cdot e
\]
defines a right logarithmic connection on $(E,B)$ such that
\[\cE\simeq \bar\cE\otimes_{\shD_{E,B}}\shD_E,\]
where $\shD_{E,B}$ is the sheaf of logarithmic differential operators along $B$ (see \cite{WZ}).
Since none of its residues along the components $E_W\cap E$ is a negative integer, $\cE$ is the $!$-extension of $\bar\cE|_{E^o}$. Then, we obtain the logarithmic Spencer resolution of $\bar\cE$:
\[\cdots\to \bigwedge^r\cT_E(-\log B)\otimes_\cO\shD_{E,B}\to \bigwedge^{r-1}\cT_E(-\log B)\otimes_\cO\shD_{E,B}\to \cdots\to \shD_{E,B}. \]
Extending scalars in the above
logarithmic Spencer resolution of $\bar \cE$ to $\shD_E$ yields a locally
free right $\shD_E$-resolution of $\cE$.
Using $e$ to identify $\bar\cE$ with $\cO_E$, we therefore
represent the derived tensor product
\[
 \cE\otimes^{\mathbb L}_{\shD_E}\cO_E
\]
by a complex $K^\bullet$ with terms
\[
 K^{-r}=\bigwedge^r\cT_E(-\log B),
 \qquad 0\le r\le n-1,
\]
whose differentials are determined by \eqref{eq:raction}
and the Leibniz rule.

Under the identification
\[
 \wt\mu_+(\cE)\simeq R\Gamma(E,K^\bullet),
\]
the image of $e$ in $\cH^0\wt\mu_+(\cE)$ is represented by
$1\in H^0(E,K^0)=H^0(E,\cO_E)$.
Consider the hypercohomology spectral sequence
\[
 E_1^{-r,j}
 =H^j\left(E,\bigwedge^r\cT_E(-\log B)\right)
 \Longrightarrow
 \mathbb H^{j-r}(E,K^\bullet).
\]
For each $\ell\ge1$, the only possible differential entering
$E_\ell^{0,0}$ is
\[
 d_\ell\colon
 E_\ell^{-\ell,\ell-1}\longrightarrow E_\ell^{0,0}.
\]
Its source vanishes: for $1\le\ell\le n-1$, Lemma
\ref{lm:NTvanishing} gives
\[
 E_1^{-\ell,\ell-1}
 =H^{\ell-1}\left(E,\bigwedge^\ell\cT_E(-\log B)\right)
 =0,
\]
while for $\ell\ge n$, the source lies outside the range
of the complex. There are also no outgoing differentials
from $E_\ell^{0,0}$, since $K^\bullet$ is concentrated in
nonpositive degrees. Thus $E_1^{0,0}=E_\infty^{0,0}$.
Moreover, the absence of positive-degree terms in $K^\bullet$
ensures that $E_\infty^{0,0}$ is a subspace of
$\mathbb H^0(E,K^\bullet)$, rather than merely a subquotient.
We therefore obtain an injection
\[
 H^0(E,\cO_E)
 =E_1^{0,0}
 =E_\infty^{0,0}
 \hookrightarrow
 \mathbb H^0(E,K^\bullet).
\]
Consequently, the image of $e$ in
$\cH^0\wt\mu_+(\cE)$ is nonzero. By Kashiwara's equivalence,
its image in $\cH^0\mu_+(\cQ_0)$ is also nonzero.
In particular,
\[
 \cH^0\mu_+(\cQ_0)\neq0.
\]

We now complete the proof of Theorem \ref{thm:stmc}.
By Lemma \ref{lm:pushlocalization}, we have an injection
\[
 \cH^0\mu_+(\widetilde{\cN}_{\mathfrak m})
 \hookrightarrow
 j_*(\cM^r_{\mathfrak m}|_U)
\]
and, consequently, an injection
\[
 \cH^0\mu_+(\cQ_0)\hookrightarrow\Psi_f^\alpha.
\]
As in the proof of Theorem \ref{thm:kashili}, under the first
injection, the image of $\tilde f^s\tilde v$ in
$\cH^0\mu_+(\widetilde{\cN}_{\mathfrak m})$ is identified with
$f^sv$. It follows that the image of
$\mu^*(f_1^{k-n})\cdot\tilde f^s\tilde v$ is identified with
$f_1^{k-n}\cdot f^sv$.

We have already shown that the image of
\[
 e=[\mu^*(f_1^{k-n})\cdot\tilde f^s\tilde v]
\]
in $\cH^0\mu_+(\cQ_0)$ is nonzero. The injection into
$\Psi_f^\alpha$ therefore implies that the image of
$f_1^{k-n}\cdot f^sv$ in $\Psi_f^\alpha$ is nonzero.
In particular, the $\shD_X[s]_{\mathfrak m}$-submodule
generated by the image of $f^sv$ is nonzero.
By Corollary \ref{cor:j!-10}, we obtain
\[
 \frac{
   f^sv\cdot\shD_X[s]_{\mathfrak m}
 }{
   f^{s+1}v\cdot\shD_X[s]_{\mathfrak m}
 }
 =
 \frac{
   f^sv\cdot\shD_X[s]_{\mathfrak m}
 }{
   j_!(\cM^r_{\mathfrak m}|_U)
 }
 \neq0.
\]
Lemma \ref{lm:z=bf} now gives $b_f(\alpha)=0$, completing
the proof of Theorem \ref{thm:stmc}.

\bibliographystyle{amsalpha}
\bibliography{mybib}
\end{document}